\documentclass[11pt]{article}
\usepackage[utf8]{inputenc}
\usepackage[a4paper, margin = 0.5in, bottom = 1in]{geometry}
\usepackage{soul}

\usepackage[utf8]{inputenc}
\usepackage[a4paper, margin = 0.5in, bottom = 1in]{geometry}
\usepackage{soul}
\usepackage{authblk}
\usepackage{minted}
\usepackage{relsize}

\date{}
\usepackage{amsmath}
\usepackage{amsfonts}
\usepackage{amssymb}
\usepackage{mathtools}
\usepackage{mathrsfs}
\usepackage{bm}
\usepackage{wasysym}
\usepackage{graphicx}
\usepackage{caption}
\usepackage{cite}
\usepackage{subcaption}
\usepackage{skak}
\usepackage[table]{xcolor}

\usepackage{amsfonts}
\usepackage{subcaption}
\usepackage{lmodern}
\usepackage{xcolor}
\usepackage{csquotes}
\usepackage{mathtools}
\usepackage{float}
\usepackage{amsmath,amssymb,amsthm,comment}

\usepackage{authblk}
\usepackage{minted}

\date{}
\usepackage{amsmath}
\usepackage{amsfonts}
\usepackage{amssymb}
\usepackage{mathtools}
\usepackage{mathrsfs}
\usepackage{bm}
\usepackage{wasysym}
\usepackage{graphicx}
\usepackage{caption}
\usepackage{cite}
\usepackage{subcaption}
\usepackage{skak}
\usepackage[table]{xcolor}

\usepackage{amsfonts}
\usepackage{subcaption}
\usepackage{lmodern}
\usepackage{xcolor}
\usepackage{csquotes}
\usepackage{mathtools}
\usepackage{float}
\usepackage{amsmath,amssymb,amsthm,comment}
\usepackage{amsfonts}
\usepackage{subcaption}
\usepackage{xcolor}
\usepackage{csquotes}
\usepackage{mathtools}
\usepackage{float}
\usepackage{amsmath,amssymb,amsthm}
\usepackage[colorlinks,citecolor=red,urlcolor=blue,linkcolor=blue]{hyperref}
\usepackage{mathrsfs,dsfont}
\usepackage{graphicx}
\usepackage{enumerate}
\usepackage{nicefrac}
\usepackage[font={normalsize},width=1.12\linewidth]{caption}
\usepackage{comment}
\usepackage{hyperref}
\usepackage{pgfplots}
\pgfplotsset{compat=1.15}
\usepackage{mathrsfs}
\usetikzlibrary{arrows}

\newtheorem{nummer}{ }[section]
\newtheorem{thm}[nummer]{Theorem}
\newtheorem{prp}[nummer]{Proposition}
\newtheorem{lem}[nummer]{Lemma}
\newtheorem{cor}[nummer]{Corollary}

\newcommand{\sx}{\mathbf S_x}
\newcommand{\sr}{\mathbf S_r}

\makeatletter
\def\opargproof[#1]{\par\noindent {\bf #1 }}
 \makeatother

\begin{document}
% \large 
% \maketitle
\medskip\medskip
\begin{center}
\vspace*{50pt}
{\LARGE\bf On the Sidon tails of $\left\{\lfloor x^n\rfloor\right\}$}

\bigskip
{\small Sayan Dutta}\\[1.2ex] 
 {\scriptsize %Département de mathématiques et de statistique, Université de Montréal\\
sayandutta345@gmail.com\\
\href{https://sites.google.com/view/sayan-dutta-homepage}{https://sites.google.com/view/sayan-dutta-homepage}}
\\[3ex]

\end{center}

\begin{abstract}
    We prove that the tail of the sets
    $$\sx := \big\{\left\lfloor x^n\right\rfloor : n\in \mathbb N\big\}$$
    are Sidon for almost all $x\in (1,2)$. Then we prove that for all $\varepsilon>0$, there exists $x\in (1,\, 1+\varepsilon)$ and $r\in (2-\varepsilon,\, 2)$ such that $\sx$ and $\sr$ do not have a Sidon tail.
\end{abstract}

\section{Introduction}
A set of positive integers $A\subset \mathbb N$ is called a \textit{Sidon Set} or a \textit{Sidon Sequence} if the equation $a+b=c+d$ does not have any non-trivial solutions in $A$. In other words, any solution $a+b=c+d$ implies $\{a,b\}=\{c,d\}$. An extensive amount of work has been done about these sets. For a short summary of recent works, see the \textit{Introduction} section of \cite{balu}. For a detailed exposition, see \cite{lit}.

It is known (and fairly easy to prove) that the sets
$$\sx := \big\{\left\lfloor x^n\right\rfloor : n\in \mathbb N\big\}$$
are Sidon for $x\ge 2$. See the remark after the proof of Theorem \ref{almost-all} for an (unnecessarily complicated) alternate proof. It seems that the question of exploring the range $1<x<2$ has not been addressed before. The main goal of the next two sections is to fill up this gap.

Of course, we have to first find a proper interpretation of this problem in the `sparse' range. We will call an $\sx$ \textit{tail Sidon} if there exists $N_0$ such that the set
$$\big\{\left\lfloor x^n\right\rfloor : n\ge N_0\big\}$$
is a Sidon set. We will also abuse notation and often shorten \textit{tail Sidon} to just \textit{Sidon}.

Our first goal is to show that for almost all $x\in (1,2)$, the set $\sx$ is Sidon. To do so, we prove an even stronger statement.

\begin{thm}\label{almost-all}
    If $\mathbf S_x$ is not Sidon, then $x$ is algebraic. In particular, one can produce integers $u> w\ge v\ge 1$ such that
    $$1+y^u-y^v-y^w=0 \qquad\mathrm{or}\qquad 1-y^v-y^w=0$$
    for $y:=\frac 1x \in (0,1)$.
\end{thm}

With this at hand, it is now imperative that we prove the existence of $x$ arbitrarily close to $1$ and $r$ arbitrarily close to $2$ for which $\sx$ and $\sr$ are not Sidon. This will be arrived at through a neat construction, appealing to Fibonacci polynomials.

\begin{thm}\label{infty}
    For all $\varepsilon>0$, there exists $x\in(1,\,1+\varepsilon)$ and $r\in (2-\varepsilon, 2)$ for which $\sx$ and $\sr$ are not Sidon.
\end{thm}
The first part essentially follows immediately from the second part, once we prove the following simple, yet crucial lemma.
\begin{lem}\label{chain}
    If $x_0>1$ is such that $\mathbf S_{x_0}$ is not Sidon, then for every integer $k\ge 2$, the number
    $$x_k:=x_0^{1/k}$$
    also satisfies that $\mathbf S_{x_k}$ is not Sidon.
\end{lem}
\begin{proof}
    For every $n\in\mathbb N$, we have
    $$\big\lfloor x_k^{kn}\big\rfloor = \left\lfloor\left(x_0^{1/k}\right)^{kn}\right\rfloor = \big\lfloor x_0^n\big\rfloor$$
    implying the subsequence $\big\{\big\lfloor x_k^{kn}\big\rfloor:n\in\mathbb N\big\}$ coincides termwise with $\big\{\big\lfloor x_0^n\big\rfloor:n\in\mathbb N\big\}$.
\end{proof}

With this at hand, (both parts of) Theorem \ref{infty} is now a corollary of the following proposition.
\begin{prp}\label{fib}
    For a positive integer $k$, let
    $$f_k(x):=\ x^k-x^{k-1}-x^{k-2}-\cdots-x-1$$
    be the $k$-Fibonacci polynomial. Let $\tilde\alpha_k$ is the unique positive real root of $f_k$. Then, $\sr$ is not Sidon with $r=\tilde\alpha_k$ for all odd $k\ge 11$.
\end{prp}

\section{$\sx$ that are Sidon}
As advertised, we will now provide a proof of Theorem \ref{almost-all}. The idea is that, for $\sx$ to not be Sidon, it would have to fall into one of two strict repeating structures.
\begin{proof}[Proof of Theorem \ref{almost-all}]
    Assume an $x$ for which $\mathbf S_x$ is not Sidon. Note that $\mathbf S_x$ is eventually increasing - so, from now on, we will only work in this increasing tail. Let $\mathbf x_n:=\left\lfloor x^n\right\rfloor$. It is easy to show that if there are indices $a$, $b$, $c$ and $d$ such that $\mathbf x_d+ \mathbf x_a = \mathbf x_b+ \mathbf x_c$, then $a<b\le c<d$.

    Now, write $\theta_n := x^n-\lfloor x^n\rfloor \in [0,1)$. If $\mathbf x_d+\mathbf x_a=\mathbf x_b+\mathbf x_c$, then
    $$x^d+x^a-x^b-x^c = (\theta_d+\theta_a)-(\theta_b+\theta_c)$$
    implying $\bigl|x^d+x^a-x^b-x^c\bigr|<2$.

    Now, define
    $$y:=\frac1x\in(0,1),\qquad u:=d-a,\quad w:=d-b,\quad v:=d-c$$
    so that
    $$\bigl|1+y^u-y^w-y^v\bigr|<2y^d$$
    for infinitely many quadruples $(u,v,w,d)$ with $d\to\infty$.

    Choose a $D_0$ for which $10y^{D_0}<1$, and from now on assume $d>D_0$. Consider any collision $-2y^d < 1+y^u-y^w-y^v < 2y^d$. This gives
    $$y^v+y^w > 1+y^u-2y^d \ge 1-2y^d > \frac12$$
    since $y^u\ge 0$ and $2y^d<\tfrac12$. Since $v\le w$ and $0<y<1$, we have $y^v\ge y^w$, hence $y^v+y^w \le 2y^v$. Therefore $2y^v>\tfrac12$, i.e. $y^v>\tfrac14$. Because $v\mapsto y^v$ is strictly decreasing, there is a bound $V=V(y)$ such that for every collision with $d\ge D_0$, one has $v\le V$.

    On the other hand, $y^v\le y$ implying
    $$y^w>1-y-2y^d>\frac 12$$
    and hence we also get a bound $W=W(y)$ such that for every collision with $d\ge D_1$, one has $w\le W$. In other words, all sufficiently large collisions lie in the finite set $\{(v,w)\in\mathbb N^2:\ 1\le v\le w\le W\}$. Since there are infinitely many collisions with $d\to\infty$, by the pigeonhole principle there exist fixed integers $v_0$, $w_0$ and an infinite subsequence of collisions along which $(v,w)=(v_0,w_0)$ and $d\to\infty$.

    Let $C:=C(y):=y^{v_0}+y^{w_0}$ so that
    $$\bigl|1+y^{u}-C\bigr|<2y^{d}\xrightarrow[]{d\to\infty}0$$
    so that $y^{u}\to C-1$.

    If $C=1$, then $y^{v_0}+y^{w_0}=1$, i.e., $1-y^{v_0}-y^{w_0}=0$.

    If $C>1$, then $y^{u}\to C-1>0$. Since $y^n\to 0$, the exponents $u$ cannot tend to $\infty$; more formally, choose $M$ such that $2y^M<C-1$. For all sufficiently large terms of the subsequence, $2y^{u}>C-1$ implying $u<M$. Thus, $u$ takes values in a finite set $\{1,2,\dots,M-1\}$, so there exists $u_0$ and a further infinite subsequence with $u=u_0$ constant.On that further subsequence, we have
    $$\bigl|1+y^{u_0}-C\bigr|<2y^{d}\xrightarrow[]{d\to\infty}0$$
    implying $1+y^{u_0}-C = 0$, and hence, $1+y^{u_0}-y^{v_0}-y^{w_0}=0$.

    This completes the proof in both cases.
\end{proof}

\textit{Remark}: Notice that the proof never really uses the fact that $x$ is in $(1,2)$. Since $y^v+y^w>1$ implies $y>1/2$, this proof also establishes that $\sx$ is Sidon for all $x>2$.

On the other hand, Theorem \ref{almost-all} and Theorem \ref{infty} meets elegantly at the interval $(\varphi,2)$. In other words, if any $x\in (\varphi,2)$ produces an $\sx$ that is not Sidon, then it has to be among the ones predicted in Proposition \ref{fib}.

\begin{cor}\label{phi}
    If $x\in (\varphi,2)$ is such that $\sx$ is not Sidon, then $x$ is in the set $\left\{\tilde\alpha_{k} : k\ge 3\right\}$.
\end{cor}
\begin{proof}
    Following Theorem \ref{almost-all}, we have $A>B\ge C\ge 1$ such that $1+y^A=y^B+y^C$ or $y^B+y^C=1$ for $y:=1/x$. But the second case is not possible because $\varphi<x<2$ implies $2y>1$, $y+y^2<1$ and $2y^2<1$.

    In the first case, since $1+y^A>1$, these inequalities also forces $B=C=1$, implying $x^{k+1}+1=2x^k$. But
    $$(x-1)f_k(x)=x^{k+1}-2x^k+1$$
    and since $x>1$, this forces $f_k(x)=0$ thus completing the proof.
\end{proof}

It should not come as a surprise to the reader that the set $(\varphi,2)$ in the statement of Corollary \ref{phi} cannot be replaced by $(1,2)$. To make this observation explicit, we will exhibit another family of polynomials in Theorem \ref{plastic} whose roots also generate Sidon collisions. More about this discussion can be found in Section \ref{con}.

\section{Existence of not Sidon $\sx$}
As advertised, we will provide a proof of Proposition \ref{fib}. For that, we need a few lemmas that are well known and often quite easy to prove (see for example \cite{sss, luca1, luca2, dub}). For the sake of completeness, we will provide proofs of these lemmas after the proof of the proposition.
\begin{lem}\label{fib-roots}
    We have $1<\tilde\alpha_k<2$, $\tilde\alpha_k^k(2-\tilde\alpha_k)=1$, and $\tilde\alpha_k\to2$ as $k\to\infty$. Moreover, $\tilde\alpha_k$ is the only root of $f_k$ outside $\{|z|<1\}$, and all roots of $f_k$ are simple. If $k$ is odd, then $\tilde\alpha_k$ is the only real root of $f_k$.
\end{lem}

\begin{lem}\label{fib-angle}
    Let $k\ge 3$ be odd.  Among the roots of $f_k$ other than $\tilde\alpha_k$, there is a conjugate pair $\beta=\rho e^{i\omega}$ and $\overline{\beta}=\rho e^{-i\omega}$ of maximal modulus, satisfying $0<\omega<\frac{5\pi}{2k}$. If $\delta\in(0,2\pi)$ satisfies $\delta\equiv k\omega\pmod{2\pi}$, then $0<\delta<\frac{\pi}{2}$.
\end{lem}

\begin{lem}\label{trace-floor}
    For $n\ge 0$, define
    $$T_n:=\sum_{j=1}^k\alpha_j^n, \qquad E_n:=\sum_{j=2}^k\alpha_j^n$$
    where $\alpha_1=\tilde\alpha_k,\alpha_2,\ldots,\alpha_k$ are the roots of $f_k$. Then, $T_n\in\mathbb Z$ for every $n\ge0$, and $T_{n+k+1}+T_n=2T_{n+k}$. Moreover, there exists $N=N(k)$ such that
    $$\left\lfloor\tilde\alpha_k^n\right\rfloor=T_n-u_n, \qquad
    u_n=
    \begin{cases}
        1,&E_n>0,\\
        0,&E_n\le0.
    \end{cases}$$
    for every $n\ge N$
\end{lem}

\begin{lem}\label{rotation}
    Let $0<\omega<2\pi$, and let $J$ be an open interval modulo $2\pi$ whose length is greater than $\omega$. Then $\{n\omega\pmod {2\pi}:n\ge1\}$ meets $J$ infinitely often.
\end{lem}

With these at hand, we can now present the promised proof.
\begin{proof}[Proof of Proposition \ref{fib}]
    Fix an odd integer $k\ge 11$, and write $\alpha:=\tilde\alpha_k$. Let $\beta=\rho e^{i\omega}$ and $\overline{\beta}$ be the conjugate pair supplied by Lemma \ref{fib-angle}. Let $\delta\equiv k\omega\pmod{2\pi}$ for $0<\delta<\frac{\pi}{2}$.

    Since $\omega<5\pi/(2k)$, we have $\pi-\delta-2\omega>\frac{\pi}{2}-\frac{5\pi}{k}>0$. Choose $\eta>0$ so small that $4\eta<\pi-\delta-2\omega$ and define
    $$J:=\left(-\frac{\pi}{2}+\eta,\ \frac{\pi}{2}-(\delta+\omega)-\eta\right)$$
    so that $|J|=\pi-\delta-\omega-2\eta>\omega$. By Lemma \ref{rotation}, there are infinitely many $n$ such that $n\omega\pmod{2\pi}\in J$. For each such $n$, we have
    \begin{align*}
        &\cos(n\omega)\ge\sin\eta\\
        &\cos((n+k)\omega)\ge\sin\eta\\
        &\cos((n+k+1)\omega)\ge\sin\eta
    \end{align*}
    since $k\omega\equiv\delta\pmod{2\pi}$.

    Let $\gamma_1,\dots,\gamma_{k-3}$ be the roots of $f_k$ other than $\alpha,\beta,\overline{\beta}$. Set $\rho_2:=\max_j|\gamma_j|<\rho$. Then
    \begin{align*}
        E_m &=\ \beta^m+\overline{\beta}^{\,m} \,+\,\sum_{j=1}^{k-3}\gamma_j^m\\
        &=2\rho^m\cos(m\omega)+R_m
    \end{align*}
    where $|R_m|\le(k-3)\rho_2^m$. Consequently, we have
    $$(k-3)\left(\frac{\rho_2}{\rho}\right)^m<\ \sin\eta$$
    for all sufficiently large $m$.

    For sufficiently large $n$ among the infinitely many indices above, this gives $E_n>0$, $E_{n+k}>0$, $E_{n+k+1}>0$. Taking $n$ larger than the integer $N$ in Lemma \ref{trace-floor}, we obtain $u_n=u_{n+k}=u_{n+k+1}=1$. Therefore, we have
    \begin{align*}
        \left\lfloor\alpha^{n+k+1}\right\rfloor +\left\lfloor\alpha^n\right\rfloor
        &=(T_{n+k+1}-1)+(T_n-1)\\
        &=2T_{n+k}-2\\
        &=2\left\lfloor\alpha^{n+k}\right\rfloor
    \end{align*}
    which is a non-trivial Sidon collision, provided $\left\lfloor\alpha^n\right\rfloor$, $\left\lfloor\alpha^{n+k}\right\rfloor$, $\left\lfloor\alpha^{n+k+1}\right\rfloor$ are distinct for sufficiently large $n$. This is guaranteed by the fact that the sequence $\left\{\lfloor\alpha^m\rfloor\right\}$ is eventually strictly increasing.
\end{proof}

\textit{Remark}: With reference to Corollary \ref{phi}, we should also note that in the interval $(\varphi,2)$, the chain of inequalities
$$\left\lfloor x^a\right\rfloor + \left\lfloor x^b\right\rfloor \le \left\lfloor x^{b-1}\right\rfloor + \left\lfloor x^b\right\rfloor \le \left\lfloor x^{b+1}\right\rfloor \le \left\lfloor x^d\right\rfloor < \left\lfloor x^d\right\rfloor + \left\lfloor x^c\right\rfloor$$
shows that the only collisions in this interval are of the type produced in the proof of Proposition \ref{fib}.

Finally, we can express the proof of the main theorem of this section.
\begin{proof}[Proof of Theorem \ref{infty}]
    The first part follows from Lemma \ref{chain} and the case $k=314159$ of Proposition \ref{fib}. The second part is true since Lemma \ref{fib-roots} guarantees $\tilde\alpha_k\to2$.
\end{proof}

Finally, to complete this section, we prove the lemmas that we used in the proof of Proposition \ref{fib}.

\begin{proof}[Proof of Lemma \ref{fib-roots}]
    Set $q_k(z):=(z-1)f_k(z)=z^{k+1}-2z^k+1$. Note that
    $$f_k(x)=\frac{1-x^k(2-x)}{x-1}$$
    for $x>0$, $x\ne1$. Also, $f_k(x)<0$ for $0<x\le 1$.

    For $1<x<2$, define $h(x):=x^k(2-x)$, implying $h'(x)=x^{k-1}\bigl(2k-(k+1)x\bigr)$. Thus $h$ increases on $\left(1,\frac{2k}{k+1}\right)$ and decreases on $\left(\frac{2k}{k+1},2\right)$. Since $h(1)=1$, $h\left(\frac{2k}{k+1}\right)>1$, and $h(2)=0$, the equation $h(x)=1$ has exactly one solution in $(1,2)$.  This is $\tilde\alpha_k$.  For $x\ge2$, we have $q_k(x)>0$, so there are no further positive roots.

    On the other hand, we also have
    $$f_k\left(\frac32\right)=2-\left(\frac 32\right)^k<0$$
    implying $\tilde\alpha_k>\frac 32$. Also, $q_k(\tilde\alpha_k)=0$ implies $\tilde\alpha_k^k(2-\tilde\alpha_k)=1$ and hence
    $$0<2-\tilde\alpha_k=\tilde\alpha_k^{-k}<\left(\frac 23\right)^k$$
    implying $\tilde\alpha_k\to2$.

    Now fix $1<r<\tilde\alpha_k$. Since $f_k(r)<0$, we have $r^{k+1}+1<2r^k$. This implies
    $$\left|z^{k+1}+1\right| \le\ r^{k+1}+1 <\ 2r^k \,=\, \left|2z^k\right|$$
    on $|z|=r$. By Rouch\'e's theorem, $q_k(z)=z^{k+1}-2z^k+1$ has exactly $k$ zeros in $|z|<r$. Since $\tilde\alpha_k>r$ is a zero, it follows that $\tilde\alpha_k$ is the only zero of $q_k$ outside $|z|<r$. Letting $r\downarrow 1$, every zero of $q_k$ other than $\tilde\alpha_k$ has modulus at most $1$.

    If $|z|=1$ and $q_k(z)=0$, then $2z^k=z^{k+1}+1$. Equality holds in the triangle inequality on the right, so $z^{k+1}=1$. The displayed equation then gives $z^k=1$, and hence $z=1$. Since $f_k(1)=1-k\ne 0$, every root of $f_k$ other than $\tilde\alpha_k$ lies strictly inside the unit disk.

    To prove simplicity, suppose that $z$ were a multiple root of $f_k$. Then it would also be a multiple root of $q_k$, and hence $q_k'(z)=z^{k-1}\bigl((k+1)z-2k\bigr)=0$. Since $f_k(0)\ne 0$, this would imply $z=\frac{2k}{k+1}$.  But $h$ is strictly larger than $1$ at its maximum $\frac{2k}{k+1}$. So, we have
    $$q_k\left(\frac{2k}{k+1}\right) = 1-h\left(\frac{2k}{k+1}\right) <0$$
    which is a contradiction.

    Finally, if $k$ is odd and $t>0$, then
    $$f_k(-t) \,=\, -t^k-\sum_{j=0}^{k-1}(-t)^j \,=\, -t^k-\frac{1+t^k}{1+t}\ <0$$
    implying $f_k$ has no negative real root, and $\tilde\alpha_k$ is its only real root.
\end{proof}

\begin{proof}[Proof of Lemma \ref{fib-angle}]
    For $0<\theta<\pi$ and $0\le s\le1$, define
    $$F(s,\theta):=s^k\left|2-se^{i\theta}\right|$$
    so that $F(s,\theta)^2=s^{2k}\bigl(4-4s\cos\theta+s^2\bigr)$. This implies
    $$\frac{\partial}{\partial s}F(s,\theta)^2 = 2s^{2k-1}\bigl(4k-(4k+2)s\cos\theta+(k+1)s^2\bigr) \,>\, 0$$
    since $\cos\theta\le 1$ and $0<s\le 1$. Thus $F(s,\theta)$ is strictly increasing in $s$.
    
    Also, $F(0,\theta)=0$ and $F(1,\theta)=\left|2-e^{i\theta}\right|>1$. Consequently, for every $\theta\in(0,\pi)$, there is a unique $s(\theta)\in(0,1)$ such that $F(s(\theta),\theta)=1$.

    On the other hand, we have
    $$\frac{\partial}{\partial\theta}F(s,\theta)^2 \,=\, 4s^{2k+1}\sin\theta \,>\, 0$$
    and so $F(s,\theta)$ is strictly increasing in $\theta$. It follows that $s(\theta)$ is strictly decreasing in $\theta$.

    Now, define
    $$\phi(\theta):=\arg\left(2-s(\theta)e^{i\theta}\right)$$
    where $\arg$ refers to the principal argument. Since the real part of $2-s(\theta)e^{i\theta}$ is positive, $\phi$ is continuous with $-\frac{\pi}{2}<\phi(\theta)<0$ for $0<\theta<\pi$. Extend continuously by setting $\phi(0)=\phi(\pi)=0$.

    Define $H(\theta):=k\theta+\phi(\theta)$ so that $H(0)=0$ and $H(\pi)=k\pi$. Since $k\ge3$, the intermediate value theorem gives some $\theta_0\in(0,\pi)$ with $H(\theta_0)=2\pi$.

    For $z=s(\theta_0)e^{i\theta_0}$, we have $\left|z^k(2-z)\right|=1$ and $\arg\bigl(z^k(2-z)\bigr)=2\pi$. Therefore $z^k(2-z)=1$ implying $q_k(z)=0$. Since $z\ne 1$, we have $f_k(z)=0$. Now let $z=se^{i\theta}$ be any root of $f_k$ in the upper half-plane. By Lemma \ref{fib-roots}, $s<1$, and $z^k(2-z)=1$.

    Write $\phi:=\arg(2-se^{i\theta})\in\left(-\frac{\pi}{2},0\right)$. We claim that $-\phi<\theta$. If $\theta\ge \frac{\pi}{2}$, this is immediate. If $0<\theta<\frac{\pi}{2}$, then
    $$\tan(-\phi) \,=\, \frac{s\sin\theta}{2-s\cos\theta} \,<\, \frac{\sin\theta}{\cos\theta} \,=\, \tan\theta$$
    because $s\cos\theta<1$. So, $-\phi<\theta$. So, $k\theta+\phi=2\pi h$ for some integer $h\ge1$. If $h=1$, then $\frac{2\pi}{k}<\theta<\frac{5\pi}{2k}$. If $h\ge 2$, then $\theta>\frac{4\pi}{k}$.
    
    Therefore the roots corresponding to $h=1$ have strictly smaller arguments than all roots corresponding to $h\ge 2$. Since $s(\theta)$ is strictly decreasing, the roots of maximal modulus are precisely a conjugate pair arising from $h=1$. Write this pair as $\beta=\rho e^{i\omega}$, $\overline{\beta}=\rho e^{-i\omega}$. Then, $0<\omega<\frac{5\pi}{2k}$.

    Finally, for $\beta$ we have $k\omega+\arg(2-\beta)=2\pi$. Since $2-\beta$ lies in the fourth quadrant, $-\frac{\pi}{2}<\arg(2-\beta)<0$. This implies $\delta:=k\omega-2\pi=-\arg(2-\beta)$ satisfies $0<\delta<\frac{\pi}{2}$.
\end{proof}

\begin{proof}[Proof of Lemma \ref{trace-floor}]
    We have
    $$T_n=\operatorname{tr}\left(C_k^n\right)\in\mathbb Z$$
    since the companion matrix of $f_k$ is an integer matrix whose eigenvalues are $\alpha_1,\dots,\alpha_k$.

    Every root $z$ of $f_k$ satisfies $z^{k+1}-2z^k+1=0$. Multiplying by $z^n$ and summing over all roots gives $T_{n+k+1}+T_n=2T_{n+k}$.

    By Lemma \ref{fib-roots}, $|\alpha_j|<1$ for $j\ge 2$. Hence $E_n\to 0$. Choose $N$ so that $|E_n|<1/10$ for all $n\ge N$. Since $\tilde\alpha_k^n=T_n-E_n$, if $E_n>0$, then $T_n-1<\tilde\alpha_k^n<T_n$, and if $E_n\le 0$, then $T_n\le\tilde\alpha_k^n<T_n+1$.
\end{proof}

\begin{proof}[Proof of Lemma \ref{rotation}]
    Assume $\omega/2\pi$ is irrational (since otherwise, the proof is trivial). For every $M\ge1$, the pigeonhole principle gives $1\le q\le M$ such that
    $$0<\left\|\frac{q\omega}{2\pi}\right\|<\frac 1M$$
    implying the orbit is dense modulo $2\pi$, since $M$ is arbitrary. Every tail of the orbit is a translate of the whole orbit and is therefore also dense, so every nonempty open interval is visited infinitely often.
\end{proof}

% \section{Conclusion}\label{con}
% The natural question to ask here is whether it is possible to find a clean characterization of every $x\in (1,2)$ for which $\sx$ is Sidon. From personal calculations, it is possible to make the following conjecture.
% \begin{conj}
%     Let $\mathcal S$ be the set of all $x\in(1,2)$ for which $\sx$ is Sidon. Define
%     \begin{align*}
%         &B_{u,v,w}(x):=x^u-x^v-x^w+1\\
%         &G_{u,v,w}(x) := \frac{B_{u,v,w}(x)}{x-1} = x^v+x^{v+1}+\cdots+x^{u-1} -\left(1+x+\cdots+x^{w-1}\right)
%     \end{align*}
%     for integers $u>v\ge w\ge1$ and $u<v+w$. Define
%     $$\mathcal M:=\tilde{\mathcal M}\setminus \{x^2-x-1\}$$
%     where $\tilde{\mathcal M}$ is the set of all monic irreducible polynomials $m\in\mathbb Z[x]$ such that $m(x)\mid G_{u,v,w}(x)$ for at least one admissible triple $u,v,w$. For an algebraic number $x$, let $m_x$ denote its monic polynomial over $\mathbb Q$. Then, $x\notin \mathcal S$ if and only if $x$ is algebraic or $m_x\in \mathcal M$.
% \end{conj}

\section*{Acknowledgement}
A previous version of the paper had much weaker replacements of Theorems \ref{almost-all}, \ref{infty} and \ref{plastic}. ChatGPT 6 Astra found the stronger versions presented in this note, often with simpler arguments than previously appealed to. The original note can be found by checking out the older versions of this paper on arXiv.

I want to express my heartfelt thanks to my supervisor, Prof. Ramachandran Balasubramanian, who introduced me to the beautiful topic of Sidon sets. His straightforward approach to understanding mathematics has greatly influenced my own mathematical intuition and helped me greatly in writing this note.

I would also like to thank the referee for their useful suggestions.

\section*{Appendix A : Another Family of Sidon collisions}\label{app1}
To continue the discussion after Corollary \ref{phi}, we will provide another family of polynomials that also generate Sidon sets. As opposed to Proposition \ref{fib}, here we actually see non-trivial collisions of the type $a+b=c+d$ with distinct $a$, $b$, $c$, $d$. The reader is urged to observe the similarity in the overall structure of this proof and the proof of Proposition \ref{fib}.

\begin{thm}\label{plastic}
    For a positive integer $d$, let
    $$p_d(x):=\ x^d-x^{d-2}-x^{d-3}-\cdots-x-1$$
    be the $k$-plastic polynomial. Let $\alpha_d$ be the unique positive real root of $p_d$. Then, $\sr$ is not Sidon with $r=\alpha_d$ for all $d\ge 3$.
\end{thm}

The $d=3$ case of this family admits the (so called) plastic constant as one of its roots. As was the case with the previous proof, we first need the following descriptions of $\alpha_d$, before we can jump to a proof. This family of polynomials have also been extensively studied before, and the ideas going into the proof of the following observations might very well have appeared in literature before (see for example, \cite{plastic1, plastic2, plastic3}).

\begin{lem}\label{1}
    The polynomial
    $$m_d(x):=\begin{cases}p_d(x),& d\text{ odd},\\ \displaystyle\frac{p_d(x)}{x+1},& d\text{ even}\end{cases}$$
    is the minimal polynomial of $\alpha_d$, and $\alpha_d$ is a Pisot number. Moreover, $\alpha_d$ is the only real root of $m_d$.
\end{lem}

\begin{lem}\label{2}
    Among the roots of $m_d$ other than $\alpha_d$, there is a unique pair of maximal modulus, and it has the form $\beta=\rho e^{i\omega}$, $\overline\beta=\rho e^{-i\omega}$ for $0<\omega<\pi$. Furthermore, $\omega/\pi$ is irrational.
\end{lem}

\begin{lem}\label{3}
    Let $D$ be the degree of $m_d$ and let $\lambda_1=\alpha_d,\lambda_2,\ldots,\lambda_D$ be the roots of $m_d$. Define
    $$T_n=\sum_{j=1}^D\lambda_j^n,\qquad E_n=T_n-\alpha_d^n$$
    for $n\ge 0$. Then, $T_n\in \mathbb Z$, $E_n\to 0$ and $T_{n+d+1}+T_n=T_{n+d}+T_{n+d-1}$ for all $n\ge 0$. Also, we have
    $$\lfloor\alpha_d^n\rfloor=T_n-u_n,\qquad u_n= \begin{cases}
                                                    1,& E_n>0,\\
                                                    0,& E_n\le 0
                                                    \end{cases}$$
    for all sufficiently large $n$.
\end{lem}
The fact that $T_n\in \mathbb Z$ follows directly from using Newton's identities, but we will see an explicit proof.

\begin{lem}\label{4}
    Let $\chi(t):=\mathbf 1_{\{\cos t>0\}}$. There is a nonempty open interval $J$, modulo $2\pi$, such that
    $$\chi(t)+\chi(t+(d+1)\omega) \,=\, \chi(t+d\omega)+\chi(t+(d-1)\omega)$$
    for $t\in J$, and $n\omega\pmod{2\pi}$ is in $J$ for infinitely many $n$. 
\end{lem}
The proof will establish that this fact is true for any $\omega$ that is an irrational multiple of $\pi$, and therefore for our $\omega$ (using Lemma \ref{2}).

With these observations at hand, the proof of Theorem \ref{plastic} is almost immediate.
\begin{proof}[Proof of Theorem \ref{plastic}]
    Write $\alpha=\alpha_d$, and follow the notation of Lemma \ref{3}. By Lemma \ref{2}, after reordering the roots we may take $\lambda_2=\beta=\rho e^{i\omega}$ and $\lambda_3=\overline\beta=\rho e^{-i\omega}$ where $\omega/\pi$ is irrational and all remaining roots have modulus strictly less than $\rho$. So, we have
    $$E_m=2\rho^m\cos(m\omega)+R_m, \qquad R_m=\sum_{j=4}^D\lambda_j^m$$
    so that $R_m=0$ if $D=3$. Otherwise, $|R_m|\le(D-3)\rho_2^m$ where $\rho_2=\max_{4\le j\le D}|\lambda_j|<\rho$.

    Take the interval $J$ from Lemma \ref{4}. Its closure can be chosen inside an arc where the four shifted cosine functions do not vanish. Hence there is a constant $c>0$ such that $|\cos(t+a\omega)|\ge c$ for $t\in\overline J$ and $a\in\{0,d-1,d,d+1\}$. For all sufficiently large $m$, the remainder $R_m$ is smaller in absolute value than $2c\rho^m$. Consequently, whenever $n\omega\pmod{2\pi}\in J$ and $n$ is sufficiently large, $E_{n+a}$ has the same sign as $\cos((n+a)\omega)$ for each of the four offsets $a$. We can also take $n$ large enough that $|E_{n+a}|<1$. So, by Lemma \ref{3}, we have $u_{n+a}=\chi((n+a)\omega)$ for $a=0,d-1,d,d+1$. Lemma \ref{4} therefore gives
    $$u_n+u_{n+d+1}=u_{n+d}+u_{n+d-1}$$
    for infinitely many $n$.

    For each such sufficiently large $n$, we have
    \begin{align*}
        \lfloor\alpha^{n+d+1}\rfloor+\lfloor\alpha^n\rfloor
                &=T_{n+d+1}+T_n-(u_{n+d+1}+u_n)\\
                &=T_{n+d}+T_{n+d-1}-(u_{n+d}+u_{n+d-1})\\
                &=\lfloor\alpha^{n+d}\rfloor+\lfloor\alpha^{n+d-1}\rfloor
    \end{align*}
    using Lemma \ref{3}.

    As before, distinctness of these terms follows from the fact that the sequence at hand is eventually strictly increasing.
\end{proof}

\begin{proof}[Proof of Lemma \ref{1}]
    Notice that for $x>0$, the polynomial
    $$p_d(x)=x^d\left(1-\sum_{\ell=2}^{d}x^{-\ell}\right)$$
    has exactly one real root since the summation inside is continuous and strictly decreasing. Also, $p_d(1)=2-d<0$ and $p_d(\varphi)=\varphi>0$. So, $1<\alpha_d<\varphi<2$. Also, $p_d'(\alpha_d)>0$ implying $\alpha_d$ is a simple root.

    On the other hand, it is not difficult to see that $p_d$ does not have any negative real roots for odd $d$, and $(-1)$ is the only negative real root for even $d$.

    We next locate the nonreal roots. For $0\le s<1$, the polynomial
    $$q_s(z)=z^{d-1}(z^2-z-1)+s$$
    has the same number of zeros inside the unit circle as $q_0(z)=z^{d-1}(z^2-z-1)$, namely $d$. This follows from Rouché’s theorem since on $|z|=1$, we have $|z^2-z-1|^2=5-4(\operatorname{Re}z)^2\ge 1$ implying $|s|<|z^{d-1}(z^2-z-1)|$. Since $q_s$ has degree $d+1$ and no roots on the unit circle, it has exactly one root outside the unit circle. At $s=1$, $q_1=(z-1)p_d(z)$ has at most one root outside the unit circle, counting multiplicity. Since $\alpha_d>1$ is a root of $q_1$, it is the only such root. It is now clear from the definition of $m_d$ that every root of $m_d$ other than $\alpha_d$ has modulus strictly less than $1$.

    To prove irreducibility, let $h$ be the minimal polynomial of $\alpha_d$. Since $\alpha_d$ is a root of the monic integer polynomial $m_d$, Gauss’s lemma gives $m_d=h g$ with $h,g$ monic integer polynomials. If $g$ were nonconstant, all its roots would be roots of $m_d$ other than $\alpha_d$, since $\alpha_d$ is simple. They would therefore all have modulus less than $1$. But $g(0)$ is a nonzero integer, whereas the product formula for the constant term gives $|g(0)|<1$, a contradiction. Thus $g=1$, so $m_d=h$. Its roots other than $\alpha_d$ lie strictly inside the unit circle, proving that $\alpha_d$ is Pisot.
\end{proof}

The proof of Lemma \ref{2} requires a standard fact that originally appeared as Problem 5931 in the American Mathematical Monthly, proposed by Smyth \cite{Smyth1975}, and later restated (and generalised) by Zaïmi \cite{Zaimi2016}. Fortunately, it is not too difficult to prove, and therefore it is best to include it here for the sake of completeness.
\begin{lem}\label{smyth}
    Let $\gamma>1$ be an algebraic integer of degree $n$, with conjugates $\gamma_1=\gamma,\gamma_2,\ldots,\gamma_n$, and suppose $|\gamma_i|<1$ for $2\le i\le n$. If $i\ne j$ and $|\gamma_i|=|\gamma_j|$, then $\gamma_i=\overline{\gamma_j}$.
\end{lem}

\begin{proof}[Proof of Lemma \ref{2}]
    By Lemma \ref{1}, $m_d$ has degree at least $3$, and its only real root is $\alpha_d$. Thus among the other roots there is a nonreal conjugate pair. Choose a root $\beta$ of maximal modulus $\rho$ in the upper half-plane, and write $\beta=\rho e^{i\omega}$, with $0<\omega<\pi$. Lemma \ref{smyth} implies that its conjugate $\overline\beta$ is the only other root of modulus $\rho$.

    Suppose $\omega/\pi$ were rational. Then $\zeta=\beta/\overline\beta=e^{2i\omega}$ would be a root of unity. Since $m_d$ is irreducible, an automorphism $\sigma$ of its splitting field can be chosen with $\sigma(\beta)=\alpha_d$. Put $\gamma=\sigma(\overline\beta)$. Injectivity of $\sigma$ gives $\gamma\ne\alpha_d$, so Lemma \ref{1} implies $|\gamma|<1$. But $\sigma(\zeta)=\alpha_d/\gamma$ would still be a root of unity and so have modulus $1$, whereas $|\alpha_d/\gamma|>1$. This is a contradiction.
\end{proof}

\begin{proof}[Proof of Lemma \ref{3}]
    Write
    $$m_d(x)=x^D+c_1x^{D-1}+\cdots+c_D$$
    so that
    $$T_n+c_1T_{n-1}+\cdots+c_{n-1}T_1+n c_n=0$$
    for $1\le n\le D$. For $n>D$, summing $m_d(\lambda_j)=0$ after multiplying by $\lambda_j^{n-D}$ gives $T_n+c_1T_{n-1}+\cdots+c_DT_{n-D}=0$. These identities prove integrality by induction.

    Now, all roots other than $\alpha_d$ have modulus less than $1$, so $E_n\to 0$. Also, every root $\lambda_j$ of $m_d$ is a root of $p_d$. Multiplying this identity at $\lambda_j$ by $\lambda_j^n$ and summing over $j$ gives the first required identity.

    Finally, $\alpha_d^n=T_n-E_n$. For all sufficiently large $n$, $|E_n|<1$. Since $T_n$ is an integer, taking floor gives the second one.
\end{proof}

\begin{proof}[Proof of Lemma \ref{4}]
    Consider the $2\pi$-periodic step function
    $$\Psi(t)\ =\ \chi(t) \,+\, \chi(t+(d+1)\omega) \,-\, \chi(t+d\omega) \,-\, \chi(t+(d-1)\omega)$$
    and recall that $\omega/\pi$ is irrational (by Lemma \ref{2}). The four offsets $0,d-1,d,d+1$ are distinct, and the zero sets of the four shifted cosine functions are pairwise disjoint. So, a common zero for shifts $a$ and $b$ would imply $(a-b)\omega\in\pi\mathbb Z$. Therefore $\Psi$ has distinct jump points, and at each one its value changes by $1$ or $-1$.

    Over a period, $\Psi$ has an integral of $0$, since each translate of $\chi$ has an integral $\pi$. On the open arcs between its jump points, $\Psi$ is integer-valued and constant. If none of those values were zero, the unit jumps would prevent a change from positive values to negative values; all values would have the same sign, contradicting the zero integral. Hence $\Psi=0$ on some nonempty open arc. Choose $J$ inside that arc. Then, the required property holds throughout $J$.

    Finally, $\omega/2\pi$ is irrational (by Lemma \ref{2}), so its multiples are dense modulo $1$, and every tail of that orbit is dense. Thus the orbit $n\omega\pmod{2\pi}$ visits $J$ infinitely often.
\end{proof}

\begin{proof}[Proof of Lemma \ref{smyth}]
    Let $\gamma_1=\gamma,\gamma_2,\ldots,\gamma_n$ be the conjugates, and set
    $$F(X)=\prod_{1\le r\le s\le n} \big(X-\gamma_r\gamma_s\big)\ \in\mathbb Z[X]$$
    since its coefficients are symmetric polynomials in the conjugates. Suppose
    $i\ne j$, $|\gamma_i|=|\gamma_j|$, but $\gamma_i\ne\overline{\gamma_j}$. Then
    $$\phi:=\gamma_i\overline{\gamma_i}=\gamma_j\overline{\gamma_j}\ \in(0,1)$$
    is a double root of $F$. Indeed, the two distinct unordered pairs $\{\gamma_i, \overline{\gamma_i}\}$ and $\{\gamma_j, \overline{\gamma_j}\}$ both occur in its product. As an algebraic integer strictly between $0$ and $1$, $\phi$ has a conjugate $\phi'$ with $|\phi'|>1$ (otherwise its nonzero integer norm would have modulus less than $1$). Since $F\in\mathbb Q[X]$ and $\phi$ is a double root, $\phi'$ is also a double root. But every root of $F$ outside the unit disk must be $\gamma_1\gamma_r$ for some $r$, and these products are distinct implying such roots are simple, a contradiction.
\end{proof}

\section*{Appendix B : A Note on Problem \#198}
We will conclude this discussion by slightly diverting from our original discussion. This diversion is justified (somewhat involuntarily) by the fact that it provides a strong answer to a question of Erdős, and it comes from considering similar power floor sequences that we have been considering so far.

An old question of Erdős and Graham \cite{eg, eg2} asks whether a set $A$ being Sidon forces the complement of $A$ to contain an infinite arithmetic progression. This is listed as Problem \#198 in Bloom's database of Erdős Problems \cite{bloom}. As discussed in this post, this question was answered in the negative in \cite{baum}, and then again by AlphaProof which realized that the explicit construction
$$A = \big\{ (n+1)!+n :\, n\geq 0\big\}$$
works.

Here we will prove a stronger statement. Fix integers $d\geq 1$ and $r\in\{0,1,\ldots,d-1\}$.  For $n\geq 1$, define
$$A_n(d,r):= \bigcup_{\substack{m\geq 1\\m\equiv r\pmod d}} \left(m^{1/n},(m+1)^{1/n}\right) \subset(1,\infty)$$
and put
$$U_{d,r,N}:=\bigcup_{n\geq N}A_n(d,r)$$
for $N\geq 1$. Now, if $x\in A_n(d,r)$, then for some integer $m\geq 1$ with
$m\equiv r\pmod d$ we have $m<x^n<m+1$, and hence $\lfloor x^n\rfloor=m\equiv r\pmod d$.

Each $A_n(d,r)$ is open, and therefore each $U_{d,r,N}$ is open. We now prove that $U_{d,r,N}$ is dense. Indeed, let $(a,b)\subset(1,\infty)$ be a nonempty open interval. Since $1<a<b$, we have
$$b^n-a^n = b^n\left(1-\left(\frac ab\right)^n\right) \longrightarrow\infty$$
as $n\to \infty$. So, we may choose $n\geq N$ so large that $b^n-a^n>d+1$. Let $m$ be the smallest integer strictly larger than $a^n$ which is congruent to $r$ modulo $d$. So $a^n<m\leq a^n+d<b^n-1$ implying $a<m^{1/n}<(m+1)^{1/n}<b$. This gives
$$\bigl(m^{1/n},(m+1)^{1/n}\bigr)\subset(a,b)$$
and this interval is contained in $A_n(d,r)\subset U_{d,r,N}$. Since $(a,b)$ was arbitrary, $U_{d,r,N}$ is dense in $(1,\infty)$.

It follows from the Baire Category Theorem \cite{fol} that
$$\mathcal G_{d,r}:=\bigcap_{N=1}^{\infty}U_{d,r,N}$$
is a residual subset of $(1,\infty)$. Since there are only countably many pairs
$(d,r)$, we have that
$$\mathcal G:= \bigcap_{d=1}^{\infty}\ \bigcap_{r=0}^{d-1}\mathcal G_{d,r}$$
is residual in $(1,\infty)$.

Let $x\in\mathcal G$. Fix $d\geq 1$ and $r\in\{0,\ldots,d-1\}$. Since $x\in U_{d,r,N}$ for every $N$, for every $N$ there exist an exponent $n\geq N$ and an integer $m\geq 1$ such that $\lfloor x^n\rfloor=m\equiv r\pmod d$. Choosing these exponents successively, we obtain a subsequence $\lfloor x^{n_j}\rfloor\to\infty$. Hence $\sx$ contains infinitely many, and indeed arbitrarily large, integers congruent to $r$ modulo $d$. This proves the first assertion.

Now, let $P=\{a+kd:k\geq 0\}$ be an infinite arithmetic progression, where $a\geq 1$ and $d\geq 1$. Taking $r\equiv a\pmod d$, the preceding assertion gives an element
$m\in\sx$ with $m\equiv a\pmod d$ and $m\geq a$. Thus $m=a+kd$ for some $k\geq 0$, so $m\in P$. Therefore $\sx$ meets every infinite arithmetic progression, and $\mathbb N\setminus\sx$ contains none.

Consequently, every $x\in\mathcal G\cap(2,\infty)$ produces a Sidon set $\mathbf S_x$ whose complement contains no infinite arithmetic progression. Since $\mathcal G\cap(2,\infty)$ is residual in $(2,\infty)$, this gives a residual family of such examples.

\bibliographystyle{plain}
\bibliography{diss}

@article{balu,
title = {{The $m$-th element of a Sidon set}},
journal = {Journal of Number Theory},
volume = {279},
pages = {594-602},
year = {2026},
issn = {0022-314X},
note = {\url{https://doi.org/10.1016/j.jnt.2025.07.007}},
author = {R. Balasubramanian and Sayan Dutta}
}

@article{lit,
  title={A complete annotated bibliography of work related to {S}idon sequences},
  author={O'Bryant, Kevin},
  journal={Electronic Journal of Combinatorics},
  volume={DS11},
  number={39},
  year={2004},
  note = {\url{https://doi.org/10.37236/32}}
}

@misc{bloom,
    author = {T. F. Bloom},
    title = {Erdős Problem \#198},
    note = {\url{https://www.erdosproblems.com/198}}
}

@misc{sss,
    author = {Saha, Satvik},
    title ={{$K$-bonacci Numbers}},
    note = {\url{https://sahasatvik.github.io/kbonacci/}}
}

@article{baum,
title = {Partitioning vector spaces},
journal = {Journal of Combinatorial Theory, Series A},
volume = {18},
number = {2},
pages = {231-233},
year = {1975},
issn = {0097-3165},
note = {\url{https://www.sciencedirect.com/science/article/pii/0097316575900163}},
author = {James E Baumgartner}
}

@article{eg,
  title={{Old and new problems and results in combinatorial number theory: van der Waerden's theorem and related topics}},
  author={Erdős, Paul and Graham, RL},
  journal={Enseign. Math.},
  pages={325-344},
  year={1979}
}

@article{eg2,
  title={Old and new problems and results in combinatorial number theory},
  author={Erdős, Paul and Graham, RL},
  journal={Monographies de L'Enseignement Mathematique},
  year={1980}
}

@book{kn,
  title={Uniform distribution of sequences},
  author={Kuipers, Lauwerens and Niederreiter, Harald},
  year={2012},
  publisher={Courier Corporation},
  note={\url{https://web.maths.unsw.edu.au/~josefdick/preprints/KuipersNied_book.pdf}}
}

@book{fol,
  title={Real Analysis: Modern Techniques and Their Applications},
  author={Folland, Gerald B},
  year={1999},
  note={\url{https://apachepersonal.miun.se/~andrli/Bok.pdf}},
  publisher={John Wiley \& Sons}
}

@article{dub,
  title={{No two non-real conjugates of a Pisot number have the same imaginary part}},
  author={Dubickas, Art{\=u}ras and Hare, Kevin and Jankauskas, Jonas},
  journal={Mathematics of computation},
  volume={86},
  number={304},
  pages={935--950},
  year={2017},
  note={\url{https://doi.org/10.1090/mcom/3103}}
}

@article{luca1,
title = {{On the separation of the roots of the generalized Fibonacci polynomial}},
author={Garc{\'\i}a, Jonathan and G{\'o}mez, Carlos A and Luca, Florian},
journal = {Indagationes Mathematicae},
volume = {35},
number = {2},
pages = {269-281},
year = {2024},
issn = {0019-3577},
note = {\url{https://doi.org/10.1016/j.indag.2023.12.002}}
}

@article{luca2,
  title={{On the arguments of the roots of the generalized Fibonacci polynomial}},
  author={Alahmadi, Adel and Klurman, Oleksiy and Luca, Florian and Shoaib, Hatoon},
  journal={Lithuanian Mathematical Journal},
  volume={63},
  number={3},
  pages={249--253},
  year={2023},
  publisher={Springer},
  note={\url{https://doi.org/10.1007/s10986-023-09604-0}}
}

@article{Smyth1975,
  title={The conjugates of algebraic integers},
  author={Smyth, C. J.},
  journal={The American Mathematical Monthly},
  volume={82},
  number={1},
  pages={86},
  year={1975},
  note={Solution to Advanced Problem 5931; \url{https://doi.org/10.2307/2319148}}
}

@article{Zaimi2016,
  title={On the conjugates of certain algebraic integers},
  author={Za\"{\i}mi, Toufik},
  journal={Publications de l'Institut Mathématique, Nouvelle Série},
  volume={99},
  number={113},
  pages={281--285},
  year={2016},
  note={\url{https://doi.org/10.2298/PIM1613281Z}}
}

@article{plastic1,
  title={The plastic number and its generalized polynomial},
  author={Iliopoulos, Vasileios},
  journal={Cogent Mathematics},
  volume={2},
  number={1},
  year={2015},
  note={Article 1023123; \url{https://doi.org/10.1080/23311835.2015.1023123}}
}

@article{plastic2,
  title={Generalized Fibonacci numbers and associated matrices},
  author={Miles, E. P., Jr.},
  journal={The American Mathematical Monthly},
  volume={67},
  number={8},
  pages={745--752},
  year={1960},
  note={\url{https://doi.org/10.1080/00029890.1960.11989593}}
}

@article{plastic3,
  title={On Generalized Fibonacci Numbers},
  author={Miller, M. D.},
  journal={The American Mathematical Monthly},
  volume={78},
  number={10},
  pages={1108--1109},
  year={1971},
  note={\url{https://doi.org/10.1080/00029890.1971.11992952}}
}

\end{document}